\documentclass[a4paper,12pt]{amsart}
\usepackage{amsmath,
			amssymb,
			amsthm,
			amsfonts,
			amscd,% Creo que no lo vamos a usar mucho, es provisional
			mathrsfs, % letras caligráficas bonitas de foliaciones
			enumerate, % para poner [i.] y todo eso
			mathtools} % Para \mathclap http://tex.stackexchange.com/questions/12344/separate-long-math-text-under-sum-symbol-into-different-lines/33358#33358
\usepackage[utf8]{inputenc}
\usepackage{psfrag}% I added them
\usepackage{graphicx}
\usepackage[all]{xy}
\usepackage[normalem]{ulem}
\usepackage{url} % para que bibtex ponga las url.
\usepackage{color}
\usepackage[top=1.5in,
			bottom=1.5in,
			left=1in,
			right=1in]{geometry}
\usepackage{braket} % para usar \set y \Set
\usepackage{etoolbox}
\makeatletter
\patchcmd{\@maketitle}
  {\ifx\@empty\@dedicatory}
  {\ifx\@empty\@date \else {\vskip3ex \centering\footnotesize\@date\par\vskip1ex}\fi
   \ifx\@empty\@dedicatory}
  {}{}
\patchcmd{\@adminfootnotes}
  {\ifx\@empty\@date\else \@footnotetext{\@setdate}\fi}
  {}{}{}
\makeatother
\newdimen\theight
\def\TeXref#1{%
             \leavevmode\vadjust{\setbox0=\hbox{{\tt
                     \quad\quad  {\small \textrm #1}}}%
             \theight=\ht0
             \advance\theight by \lineskip
             \kern -\theight \vbox to
             \theight{\rightline{\rlap{\box0}}%
             \vss}%
             }}%

\DeclareMathOperator{\im}{im}

\newtheorem{theorem}{Theorem}[section]
\newtheorem{lemma}[theorem]{Lemma}
\newtheorem{proposition}[theorem]{Proposition}

\theoremstyle{definition}
\newtheorem{definition}[theorem]{Definition}
\newtheorem{example}[theorem]{Example}

\theoremstyle{remark}
\newtheorem{remark}[theorem]{Remark}

\numberwithin{equation}{section}

\newcommand\chii{\raise2pt\hbox{$\chi$}}
\newcommand\phii{{\raise2pt\hbox{$\varphi$}}}

\newcommand\Q{\mathbb{Q}}

\newcommand\N{\mathbb{N}}
\newcommand\ese{\mathbb{S}}
\newcommand\T{\mathbb{T}}

\newcommand\Z{\mathbb{Z}}
\newcommand{\R}{\mathbb{R}}

\newcommand\F{\mathscr{F}}

\newcommand\LL{\mathscr{L}}

\newcommand\CP{\mathbb{CP}}

\newcommand{\THL}[1]{(THL)_{#1}}
\newcommand{\HL}[1]{(HL)_{#1}}
\newcommand{\HB}[1]{H_B^{#1}}
\newcommand{\HM}[1]{H_M^{#1}}
\newcommand{\PHB}[1]{P\!H_B^{#1}}
\title{Hard Lefschetz Property for Isometric Flows}

\author[J.I.~Royo Prieto]{Jos\'{e} Ignacio Royo Prieto}
\address{Matematika Saila\\ Zientzia eta Teknologia Fakultatea\\ University of the Basque Country UPV/EHU\\ Barrio Sarriena s/n\\ 48940 Leioa\\Spain.}
\email{joseignacio.royo@ehu.eus}
\thanks{The first author has been partially supported by 
	the Gobierno Vasco Grant IT1094-16; the first and second authors by Ministerio de Ciencia, Spain, grant MTM2016-77642-C2-1-P, and
	the three authors  by Ministerio de Ciencia, Spain, grant PID2019-105621GB-I00. The
	authors wish to thank the referee for the useful remarks made in order to
	improve this paper.}

\author[M.~Saralegi Aranguren]{Martintxo Saralegi-Aranguren}
\address{F\`{e}d\`{e}ration CNRS\\  Nord-Pas-de-Calais FR 2956\\ UPRES-EA 2462 LML\\
Facult\'e Jean Perrin\\ Universit\'{e} d'Artois\\   Rue Jean Souvraz SP 18\\ 62 307 Lens Cedex, France.}
\email{martin.saraleguiaranguren@univ-artois.fr}

\author[R.~Wolak]{Robert Wolak}
\address{Instytut Matematyki\\ Uniwersytet Jagiellonski\\
	ul. prof. Stanis{\l}awa {\L}ojasiewicza 6
	30-348 Krak\'ow,  Poland.}
\email{robert.wolak@im.uj.edu.pl}

\keywords{Lefschetz Hard Property, Contact Manifolds, Isometric flow}

\subjclass[2010]{53C12, 53D10, 53C25}

\date{\today}

\begin{document}

\begin{abstract}
The Hard Lefschetz Property (HLP) is an important property which has been studied in several categories of the symplectic world. 
For Sasakian manifolds, this duality is satisfied by the basic cohomology (so, it is a transverse property), but a new version of the HLP has been recently given in terms of duality of the cohomology of the manifold itself in  \cite{mino}. Both properties were proved to be equivalent (see \cite{linyi}) in the case of  K-contact flows. In this paper we extend both versions of the HLP (transverse and not) to the more general category of isometric flows, and show that they are equivalent.
We also give some explicit examples which illustrate the categories where the HLP could be considered.
\end{abstract}

\maketitle

\section*{Introduction}

The origins of the Hard Lefschetz Property (HLP in the sequel) go back to Lefschetz's  study of topological properties of algebraic real projective varieties  \cite{L24}, where he proved that the repeated cup product by the cohomology class of a hyperplane gives an isomorphism in the cohomology of the variety. Later, a version of that theorem was proved by Hodge (see \cite{hodge}) for general compact K\"ahler manifolds, stating isomorphisms between de Rham cohomology groups of complementary degrees given by multiplication by a power of the symplectic form. This property was considered to be one of the most important of this class of manifolds. Compact K\"ahler manifolds have  very strong and particular cohomological properties. A lot of effort was put into distinguishing such properties which could characterize K\"ahler manifolds within the category of compact symplectic manifolds. Now we know among other things that

\begin{itemize}
	\item there are compact symplectic manifolds which are not K\"ahler, cf. \cite{Th} for the first such an example;
	\item the torus is the only nilmanifold which is K\"ahler, cf. \cite{BG};
	\item there are compact symplectic manifolds whose cohomology ring is formal which are not K\"ahler, cf. \cite{CFL2};
	\item there are  compact Hermitian manifolds with collapsing Fr\"olicher spectral sequence which are not K\"ahler, cf. \cite{CFG2};
	\item there are compact symplectic manifolds satisfying the HLP which are not K\"ahler, cf. \cite{Y}. However, a nilmanifold having the HLP is diffeomorphic to a torus, cf. \cite{BG}, thus a K\"ahler manifold. 
\end{itemize}

Over the years many examples have been discussed and published. Among other publications let us mention the papers \cite{CFL1,CFGS,CFG1,CFG3,CFGU1,CFGU2} and the book \cite{OT}.

Within the realm of foliations the foundations of the theory of transversely K\"ahler foliations were presented  by A. El Kacimi in \cite{ELK}.  
 L. Cordero and R. Wolak in two papers presented a series of examples showing that the corresponding transverse properties  of the basic cohomology do not characterize transversely K\"ahler foliations,  cf. \cite{CW1,CW2}. 

These results proved to be of particular importance in the study of the odd-dimensional counterpart of  K\"ahler manifolds, i.e., Sasakian manifolds. In particular, by  \cite[par.~3.4.7]{ELK} the basic cohomology of a compact Sasakian manifold satisfies the HLP. In recent years a lot of research has been done to distinguish Sasakian manifolds within the class of contact metric manifolds and K-contact manifolds in particular, e.g., cf. \cite{mino,HT,MT}.

 \begin{figure}[b]
 	\psfrag{isometric}[][]{Isometric}
 	\psfrag{contact}[][]{Contact}
 	\psfrag{K}[][]{K-Contact}
 	\psfrag{lef}[][]{Lefschetz}
 	\psfrag{Sasakian}[][]{Sasakian}
 	\includegraphics[height=5.5cm]{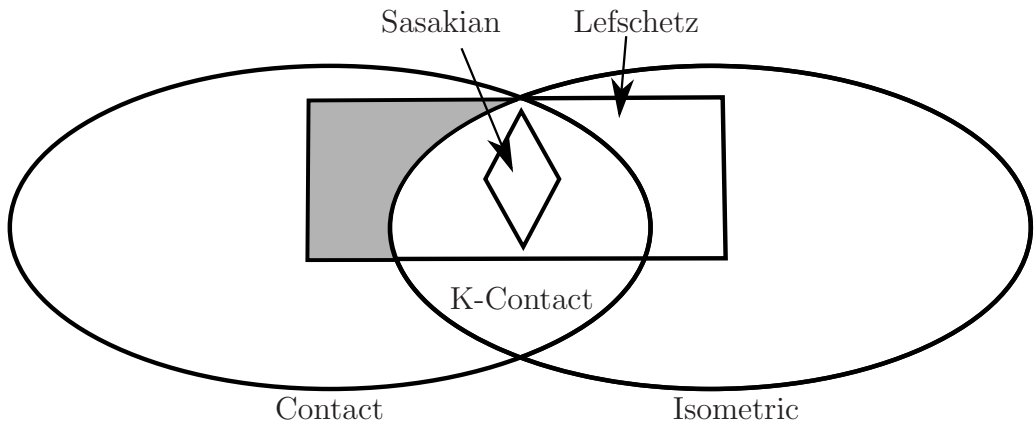}
 	\caption{Some categories where the HLP has been considered.}
 	\label{fig:dibujico}
 \end{figure}
 
One of the properties  used in these considerations was a new version of the HLP for Sasakian manifolds demonstrated in \cite{mino2} which stated Lefschetz-type isomorphisms not for the basic cohomology groups, but for the de Rham groups of the manifold itself. The authors of that paper extended the scope of the property by giving a definition of Lefschetz contact manifold in the same global terms. Examples of non-Sasakian Lefschetz contact manifolds have been given in \cite{mino2} and \cite{mino3}, all of them within the category of isometric flows (i.e. the Reeb field associated to the contact structure is a Killing vector field).

So, a priori there are two different properties (global and basic) that a contact manifold may or may not satisfy, both of them generalizing the HLP satisfied by Sasakian manifolds. In \cite{linyi} the author proves that both properties are equivalent for compact K-contact manifolds. To define the Lefschetz map the author uses the symplectic Hodge theory. We don't know whether that equivalence is held for all contact flows.

Lefschetz-type isomorphisms also exist in the realm of isometric flows, where the role of the class of the symplectic form is played by the Euler class. In this work we define two duality properties for isometric flows which resemble the HLP: a transversal one $\THL{}$  and a global one $\HL{}$. Although our definition is essentially topological and no symplectic structure is needed, in the case of K-contact flows, our new definitions agree with the previous versions of the HLP introduced above. In Section 1, we prove that both properties are equivalent for isometric flows.
So, we can call Lefschetz isometric flows the isometric flows satisfying $\THL{}$ or $\HL{}$.

In Figure \ref{fig:dibujico} we show the categories where the HLP has been defined. The HLP is satisfied in the rectangular region. We don't know whether the shaded area is nonempty (that is, whether there exist Lefschetz contact flows which are not K-contact), but all other regions are, as we illustrate with some examples in Section 2. In Example \ref{ex:cp2cp2} we provide a Lefschetz isometric flow which doesn't admit a contact structure.
In order to find an example of a flow which is contact, Lefschetz but not isometric we have to look for a flow which is not Riemannian as the Lefschetz condition ensures tautness in the Riemannian realm, thus our flow would be isometric. In the case of transversely symplectic but not Riemannian foliations  we can encounter  infinite dimensional basic cohomology which makes the Lefschetz condition problematic, as happens in Example~\ref{ex:solenoid}. We don't know whether the transversal and the global definitions of the HLP are equivalent if the contact flow is not isometric. In the referred example, neither of them is satisfied, but the problems appear at non-corresponding degrees, differently as in the isometric case.

\section{Lefschetz Duality and Transverse Duality for isometric flows}

\subsection{Preliminaries}

Throughout this section $(M,g)$ denotes a closed Riemannian manifold endowed with an isometric flow $\F$, that is, a 1-dimensional foliation defined by the orbits of a locally free $\R$-action by isometries. Let $X$ be the unit vector field defining the flow. 
$\HM{*}$ and $\HB{*}$ stand for the de Rham cohomology of $M$ and the basic cohomology of the flow, respectively. The latter is the cohomology of the complex of {\em basic forms} $\set{\omega\in\Omega(M) |  i_X\omega=i_Xd\omega=0}$. The closure of $\R$ in the group of isometries $\mathop{\rm Iso}(M,g)$ is an abelian compact and connected group, and hence, a torus $G$. We have an isomorphism $\Psi^*:H_M^*\to H^*(\Omega(M)^G)$ between the de Rham and the $G$-invariant cohomology groups (see
\cite[Th.~1 in p.~151]{GHV}).
 
We have the {\em Gysin exact sequence} (see \cite[Th.~6.13]{tondeur-libro}):
\begin{equation}\label{eq:gysink}
\xymatrix{\cdots \ar[r] 
	& \HB{k-2}\ar[r]^(.55){\varepsilon_{k-1}} 
	& \HB{k}\ar[r]^(.55){\iota_k} 
	& \HM{k}\ar[r]^(.4){\rho_k} 
	& \HB{k-1}\ar[r]^{\varepsilon_k}
	& \HB{k+1}\ar[r]
	&\cdots,}
\end{equation}
where $\iota_k$ is induced by the natural inclusion of the basic complex into the de Rham complex, $\varepsilon_k$ is the multiplication by the Euler Class $[e]=[d\chi]\in \HB{2}$ of $\F$, being $\chi=i_X g$ the characteristic form of $\F$, and $\rho_k=i_X\circ\Psi^*$, being $i_X$ the contraction operator (notice that  $\rho_k([\omega])=[i_X\omega]$ when $\omega$ is $X$-invariant). In the literature, the multiplication by the Euler class is also known as the {\em Lefschetz operator}, and is denoted by $\xymatrix@1{L:\HB{*}\ar[r]^{\wedge[e]}&\HB{*+2}}$.

\subsection{Hard Lefschetz Duality Properties}

\begin{definition}
Let $\F$ be an isometric flow on the closed manifold $M$, where $\dim{M}=2n+1$, and let $[e]\in \HB{2}$ denote its Euler class. We will say that $\F$ satisfies the {\em Transversal Hard Lefschetz property at degree $k\in\Z$} if the following property holds:
$$
\THL{k}:\quad  
L^{n-k}\colon \HB{k}\longrightarrow \HB{2n-k}\quad
\text{is an isomorphism,}
$$
where $L^{n-k}([\beta])=[\beta\wedge e^{n-k}]$. We also define the following properties:
$$
\begin{aligned}
&\THL{\le k}: &\THL{j}\text{ holds for every } j\le k\\
&\THL{}:   &\THL{j}\text{ holds for every\ } j\in\Z.
\end{aligned}
$$
In this last case, we will say that $\F$ satisfies the {\em Transversal Hard Lefschetz property}.
\end{definition}

\begin{remark}\label{rem:THLneg}
	$\THL{k}$ holds trivially if $k < 0$ or $k>2n$. For $k=0$, on one hand, if $\F$ is transversally symplectic, then $\THL{0}$ is satisfied. On the other hand, one can easily construct an $\ese^1$-principal bundle over $B=\T^4=\R^4/\Z^4$ with a nontrivial Euler class (say, $[e]=[dx_1\wedge dx_2]\in H^2_B$). We have that $[e]\ne0$, but $[e^2]=0$, and thus $\THL{0}$ does not hold.   
\end{remark}

\begin{definition}
We define the {\em $k$-th basic primitive cohomology group} as the kernel of the map $\xymatrix@1{\HB{k}\ar[rr]^{\wedge[e]^{n-k+1}}&&\HB{2n-k+2}}$, that is,
$
\PHB{k}=\Set{[\beta]\in \HB{k} |  [\beta\wedge e^{n-k+1}]=0\text{ in }\HB{2n-k+2}}.
$
\end{definition}

\begin{definition}
We will say that $\F$ satisfies the $k$-th {\em primitive condition} (and denote it by $P_k$) if the inclusion of forms induces the following two isomorphisms:
\begin{itemize}
	\item[$(P_1)_k$:]  $\xymatrix@1{i_k=\iota_k\vert_{\PHB{k}}\colon \PHB{k}\ar[r]^(0.7){\cong}& \HM{k}}$;
	\item[$(P_2)_k$:]  $\HB{k}=\PHB{k}\oplus L(\HB{k-2})$.
\end{itemize}

\end{definition}

\begin{remark}\label{remark:p0}
	Notice that as $\PHB{0}=\HB{0}=\HM{0}$, then $P_0$ is always true. $(P_1)_1$ is not always true, but for every $\beta\in\Omega^1_B$ we have $\beta\wedge e^n\in\Omega^{2n+1}_B=0$, so we have $\PHB{1}=\HB{1}$, and then $(P_2)_1$ always holds.
\end{remark}

\begin{lemma}\label{lem:iotaepi} For every $k\le n$,
$$
(THL)_{k-1}\Rightarrow \xymatrix@1{ {\iota_k}\colon \HB{k} \ar@{->}[r] & \HM{k}}\text{ is an epimorphism.}
$$
\end{lemma}

\begin{proof}
By Remark  \ref{remark:p0}, the result holds trivially for $k=0$. If $k\ge 1$, by $\THL{k-1}$, $L^{n-k+1}=L\circ L\circ\cdots\circ L$ is an isomorphism, and thus, a monomorphism. So, the first map in that composition $L=\varepsilon_k\colon \HB{k-1}\longrightarrow\HB{k+1}$ must be a monomorphism, too. The exactness of the Gysin sequence \eqref{eq:gysink} implies that $\iota_k$ is an epimorphism.
\end{proof}

\begin{remark}\label{rem:p1}
By degree reasons $\iota_1\colon\HB{1}\rightarrow\HM{1}$ is always a monomorphism (this holds for any  foliation). So, by Remark \ref{remark:p0} and Lemma \ref{lem:iotaepi}, we have that $\THL{0}$ implies $P_1$.
\end{remark}

\begin{proposition}\label{prop:pk} For every $k\le n$,
	$$
	\left.
	\begin{matrix}
	(THL)_{k-1}\\
	\text{and}\\
	(THL)_{k-2}
	\end{matrix}
		\right\}
		\Longrightarrow P_k
	$$
\end{proposition}

\begin{proof}
By Remarks \ref{remark:p0} and 	\ref{rem:p1} we have $P_0$ and $P_1$. Consider $k\ge2$. On one hand, the Gysin sequence \eqref{eq:gysink} gives
\begin{equation}\label{eq:sumadirecta}
\HB{k}\cong \im\iota_k\oplus\ker\iota_k \cong \HM{k}\oplus\im\varepsilon_{k-1}\cong \HM{k}\oplus L(\HB{k-2}),
\end{equation}
where we have used that $\iota_k$ is an epimorphism, which holds by $\THL{k-1}$ and Lemma \ref{lem:iotaepi}. On the other hand, we consider the sum $\PHB{k} + L(\HB{k-2})\le \HB{k}$. We now show that the sum is a direct one: take $[\beta]\in\PHB{k}\cap L(\HB{k-2})$. Then, there exists $[\gamma]\in\HB{k-2}$ such that 
$[\beta]=[\gamma\wedge e]\in\PHB{k}$, which implies
$$
0=[\beta\wedge e^{n-k+1}]=[\gamma\wedge e^{n-k+2}]=L^{n-k+2}([\gamma]),
$$
and by $\THL{k-2}$ we have $[\gamma]=0$. Thus, $[\beta]=0$, and the sum is direct. From \eqref{eq:sumadirecta}, we get
\begin{equation}\label{eq:oplus}
\PHB{k}\oplus  L(\HB{k-2})\le \HB{k}=\HM{k}\oplus L(\HB{k-2}),
\end{equation}
which implies that $\dim\PHB{k}\le\dim\HM{k}$. Hence, if we prove that 
$\xymatrix@1{i_k\colon \PHB{k}\ar@{->}[r]& \HM{k}}$ 
is an epimorphism, it would be an isomorphism, yielding $(P_1)_k$ and, by \eqref{eq:oplus}, $(P_2)_k$.

We complete the proof by showing that $i_k$ is an epimorphism. Let $[\alpha]\in\HM{k}$. By $\THL{k-1}$ and Lemma \ref{lem:iotaepi}, there exists $[\beta]\in\HB{k}$ such that $\iota_k([\beta])=[\alpha]$. From $\THL{k-2}$, there exists $[\gamma]\in \HB{k-2}$ such that $$
[\beta\wedge e^{n-k+1}]=L^{n-k+2}([\gamma])=[\gamma\wedge e^{n-k+2}]\in\HB{2n-k+2},
$$ 
which leads to $[(\beta-\gamma\wedge e)\wedge e^{n-k+1}]=0\in\HB{2n-k+2}$ and thus, $[\beta-\gamma\wedge e]\in\PHB{k}$. Finally, we have
$$
i_k([\beta-\gamma\wedge e])=[\beta-d(\chi\wedge\gamma)]=[\beta]\in\HM{k}.
$$
\end{proof}

\begin{definition}
Let $\F$ be an isometric flow on the closed manifold $M$, where $\dim(M)=2n+1$, and let $[e]\in \HB{2}$ denote its Euler class. We say that $\F$ satisfies the {\em Hard Lefschetz property at degree $k$} (and denote it by $(HL)_k$) if there exists an isomorphism $\LL^{n-k}\colon \HM{k}\longrightarrow \HM{2n-k+1}$ 
	making the following diagram commutative:
\begin{equation}\label{eq:diagramak}
	\xymatrix{\HM{k}\ar[rr]^(0.4){\LL^{n-k}}  
	&	& \HM{2n-k+1}\ar[d]^{\rho_{2n-k+1}}\\
		\PHB{k}\ar[u]^{i_k}\ar@{^{(}->}[r]&\HB{k}\ar[r]^{L^{n-k}}
		& \HB{2n-k}}
\end{equation}
We shall also use the following notations:
$$
\begin{aligned}
	&\HL{\le k}: &\HL{j}\text{ holds for every } j\le k\\
	&\HL{}:   &\HL{j}\text{ holds for every\ } j\in\Z.
\end{aligned}
$$
In this last case, we will say that $\F$ satisfies the {\em Hard Lefschetz property}.
\end{definition}

\begin{theorem}\label{th:equivalencia}  Let $\F$ be an isometric flow on a closed oriented manifold $M$ of dimension $2n+1$. Then, for every $k\le n$, we have $\THL{\le k} \Longleftrightarrow 
	\HL{\le k}$.
\end{theorem}
\begin{proof}	
~\\
	\fbox{$\Rightarrow$}
Assume $\THL{\le k}$. By Proposition \ref{prop:pk}, $P_k$ is true, and so $i_k$ is an isomorphism. To define $\LL^{n-k}$ we fix a basis $\set{[\beta_i]}_i$ of $\PHB{k}$. As $[\beta_i\wedge e^{n-k+1}]=0\in\HB{2n-k+2}$, there exists a basic form $\gamma_i\in\Omega_B^{2n-k+1}$ such that $\beta_i\wedge e^{n-k+1}=d\gamma_i$, and so, $\chi\wedge\beta_i\wedge e^{n-k}- \gamma_i\in\Omega^{2n-k+1}_M$ is a closed form. Thus, we can define 
\begin{equation}\label{eq:def-leches}
\LL^{n-k}([\beta_i])=[\chi\wedge\beta_i\wedge e^{n-k}- \gamma_i]
\end{equation}
and extend it by linearity. As $\gamma_i$ is basic and $\chi\wedge\beta_i\wedge e^{n-k}- \gamma_i$ is $X$-invariant, we have 
$$
\rho_{2n-k+1}([\chi\wedge\beta_i\wedge e^{n-k}- \gamma_i])=[i_X(\chi\wedge\beta_i\wedge e^{n-k}- \gamma_i)]=[\beta_i\wedge e^{n-k}],
$$ 
and the diagram  \eqref{eq:diagramak} is commutative.
Finally, from \eqref{eq:diagramak} we have 
$$
\ker\LL^{n-k}\le \ker(\rho\circ\LL^{n-k})\le\ker(L^{n-k}\vert_{\PHB{k}}\circ (i_k)^{-1})=\{0\},
$$
because $i_k$ and $L^{n-k}$ are isomorphisms. So, $\LL^{n-k}$ is a monomorphism between $\HM{k}$ and $\HM{2n-k+1}$, who have the same dimension by Poincaré Duality. Hence, an isomorphism.

~\\
\fbox{$\Leftarrow$}
First notice that $\F$ is transversally orientable, and by \cite[Th.~A]{mose85} and \cite[Th.~4.10]{eh86}, $\HB{*}$ satisfies the Poincaré Duality. In particular, $\HB{k}$ and $\HB{2n-k}$ have the same dimension. So, in order to prove that $L^{n-k}$ is an isomorphism, it suffices to show that $L^{n-k}$ is an epimorphism.

As the statement is trivial for $k<0$, we shall proceed by induction on $k$ starting at $k=-2$.
Assume that $\HL{\le k}$ holds and assume $\HL{\le k-1} \Rightarrow \THL{\le k-1}.$ 
By Proposition \ref{prop:pk}, $P_k$ holds, giving that $i_k$ is an isomorphism.
To show that $L^{n-k}$ is onto, we now take $[\varphi]\in\HB{2n-k}$. By $\THL{k-2}$,  $L^{n-k+2}$ is an isomorphism, and so, there exists $[\gamma]\in\HB{k-2}$ such that 
$$
[\varphi\wedge e]=L^{n-k+2}[\gamma]=[\gamma\wedge e^{n-k+2}]\in\HB{2n-k+2}.
$$
We have the following commutative diagram,
		\begin{equation}\label{eq:diagramaco}
		\xymatrix{\HM{k}\ar[rr]^{\LL^{n-k}}_{\cong}  
			&	& \HM{2n-k+1}\ar[d]^{\rho_{2n-k+1}}\\
			\PHB{k}\ar[u]^{i_k}_{\cong}\ar@{^{(}->}[r]&\HB{k}\ar[r]^{L^{n-k}}
			& \HB{2n-k}\ar[d]^{\varepsilon_{2n-k+1}}\\
			& \HB{k-2}\ar[r]^{L^{n-k+2}}_{\cong} & \HB{2n-k+2}}
		\end{equation}
whose right column is part of the Gysin sequence \eqref{eq:gysink}. We have
$$
[\varphi-\gamma\wedge e^{n-k+1}]\in\ker{\varepsilon_{2n-k+1}}=\im{\rho_{2n-k+1}},
$$
and so, by \eqref{eq:diagramaco} there exists $[\beta]\in\HB{k}$ such that 
$$
[\beta\wedge e^{n-k}]=L^{n-k}[\beta]=[\varphi-\gamma\wedge e^{n-k+1}],
$$
which implies
$$
[\varphi]=[\beta\wedge e^{n-k}+\gamma\wedge e^{n-k+1}]=[(\beta+\gamma\wedge e)\wedge e^{n-k}]=L^{n-k}([\beta+\gamma\wedge e]),
$$
which concludes the proof.
\end{proof}

Now, the following definition makes sense.

\begin{definition}
	We say that an isometric flow $\F$ on a closed manifold $M$ is an {\em isometric Lefschetz flow} if it satisfies $\THL{}$ or $\HL{}$.
\end{definition}

\begin{remark}
Given an isometric flow $\F$ on a compact manifold, in \cite[Section~3.2]{pitman} it is proved that the Euler classes associated to two invariant metrics are the same up to a multiplicative nonzero constant. As a result, whether an isometric flow is Lefschetz or not is a topological property in the sense that it does not depend on the chosen invariant metric, but only on the foliation $\F$ itself.
\end{remark}

In \cite{mino}, the authors define a $(2n+1)$-dimensional contact manifold $(M,\eta)$ with Reeb vector field $\xi$ to be a {\em Lefschetz contact manifold} if for every $k\le n$, the relation between $H_M^{k}$ and $H_M^{2n+1-k}$ defined by
\begin{equation}\label{eq:relation}
\mathcal{R}^k=\Set{([\beta],[\eta\wedge (d\eta)^{n-k}\beta]) | 
	\beta\in \Omega_M^k, d\beta=0, i_{\xi}\beta=0, (d\eta)^{n-k+1}\wedge\beta=0}
\end{equation}
is the graph of an isomorphism $H^k_M\cong H_M^{2n-k+1}$. 
\begin{remark}\label{rem:extremal}
Notice that $\mathcal{R}^0$ is always the graph of the isomorphism $H^0_M\cong H^{2n+1}_M$ because $(d\eta)^{n+1}=0$.
\end{remark}
We now see that if $\xi$ is Killing (i.e. we have a K-contact flow), then this notion is equivalent to $\HL{}$.

\begin{proposition}\label{prop:mino}
Let $(M,\eta)$ be a $(2n+1)$-dimensional K-contact manifold. Then, the isometric flow defined by its Reeb vector field is an isometric Lefschetz flow if and only if $(M,\eta)$ is a Lefschetz contact manifold.
\end{proposition}

\begin{proof}
We have that $X=\xi$ is a Killing vector field, $\chi=\eta$ and $e=d\eta$. If  $(M,\eta)$ is a Lefschetz contact manifold, then for every $k\le n$, the isomorphism $\LL^{n-k}$ whose graph is the relation \eqref{eq:relation} clearly makes the diagram \eqref{eq:diagramak} commutative and so, $X$ defines an isometric Lefschetz flow. Conversely, if $X$ is a Lefschetz isometric flow, it satisfies $\THL{}$ and we can construct an isomorphism $\LL^{n-k}\colon H_M^{k}\to H_M^{2n-k+1}$ as in the \fbox{$\Rightarrow$} part of Theorem \ref{th:equivalencia}. As $(M,\eta)$ is contact, by \cite[Th.~11(1)]{czarnecki}, each basic cohomology class has a harmonic representative, and thus, each primitive basic cohomology class admits a  primitive basic representative\footnote{This is \cite[Lemma~2.11]{linyi}, which can also be proved by the last paragraph of the proof of Theorem~0.1 of \cite{yan}, which applies verbatim to the complex of basic forms of the contact manifold.}. So, as $H_M^k\cong\PHB{k}$, we can find a basis   $\set{[\beta_i]}_i$ of $H_M^k$, where $\beta_i$ are primitive closed basic forms. 
In the proof of  Theorem \ref{th:equivalencia} we have to add the forms $\gamma_i$ to get  closed forms, but as  the $\beta_i$ are primitive, we can choose $\gamma_i=0$ and so  $\LL^{n-k}([\beta])=[\chi\wedge\beta\wedge e^{n-k}]$ defines an isomorphism whose graph is the relation \eqref{eq:relation}. Thus, $(M,\eta)$ is a Lefschetz contact manifold.
\end{proof}

In \cite{mino}, the authors prove that the small odd Betti numbers (up to the middle dimension) of a Lefschetz Contact flow are even. As we show now, the same algebraic proof works to prove the corresponding result for isometric Lefschetz flows.

\begin{theorem}
Let $\F$ be an isometric Lefschetz flow on the compact manifold $M$, being $\dim(M)=2n+1$. Then, the Betti number $b_k(M)$ is even for every odd $k\le n$.
\end{theorem}
\begin{proof}
Consider a basis of primitive basic classes $\set{[\beta_i]}_i$ of $H^k_M$, and construct an isomorphism $\LL^{n-k}\colon H_M^k\to H_M^{2n-k+1}$ as in the proof of Theorem~\ref{th:equivalencia}. Consider the non-degenerate bilinear form $B$ on $H_M^k$ defined as the composition
\begin{equation*}
\xymatrix@C=4pc{
	\HM{k}\otimes H_M^k
	\ar[r]^{1_M\otimes\LL^{n-k}}
	& H_M^k\otimes H_M^{2n-k+1} \ar[r]^(0.65){P}
	&\R,}
\end{equation*}
being $P([\omega_1],[\omega_2])=\int_M\omega_1\wedge\omega_2$ the usual non-degenerate pairing. Now, we have
\begin{equation*}
\begin{split}
B([\beta_i],[\beta_j])&=P([\beta_i],\LL^{n-k}[\beta_j])\\
&=\int_M(\beta_i\wedge\chi\wedge\beta_j\wedge e^{n-k} - \beta_i\wedge\gamma_i)\\
&=\int_M \beta_i\wedge\chi\wedge\beta_j\wedge e^{n-k}
\end{split}
\end{equation*}
where we have used that $\beta_i\wedge\gamma_i=0$ because it is a basic form of degree $2n+1$. As $\beta_i\chi\beta_j=(-1)^k\beta_j\chi\beta_i$, it follows that $B([\beta_i],[\beta_j])=(-1)^kB([\beta_j],[\beta_i])$. So, $B$ is a non-degenerate skew-symmetrical bilinear form, and the dimension of $\HM{k}$ must be even.
\end{proof}

\section{Examples}

\begin{example}[Sasakian manifolds]
Consider a Sasakian manifold $M$ of dimension $2n+1$ (for the essentials of Sasakian geometry we refer the reader to \cite{galicki}). Recall that its associated Reeb vector field $X$ defines an isometric flow with respect to the metric $g$ of the Sasakian structure, being the associated contact form $\chi=i_Xg$ the characteristic form of the isometric flow. As any Sasakian manifold is transversally K\"ahler, it satisfies  $\THL{}$ 
(cf.\cite[3.4.7]{ELK}), and by Theorem \ref{th:equivalencia}, it satisfies $\HL{}$, which has been proved in \cite[Section~4]{mino}.
%
%In fact, the 
%
\end{example}

In \cite{BW58} Boothby and Wang use the construction described by Kobayashi in \cite[Th.~2]{kobayashi}) to get examples of contact manifolds out of integral symplectic forms. The same construction can also be applied to get isometric flows with a prescribed integral Euler form as follows: given an integral closed form $\omega\in\Omega^2(B)$, Kobayashi's construction gives an $\ese^1$-principal bundle $\pi\colon M\to B$ whose connection form $\chi\in\Omega^1(M)^{\ese^1}$ satisfies $d\chi=\pi^*\omega$. Let $\F$ be the foliation on $M$ defined by the orbits of the principal $\ese^1$-action and consider on $TM=T\F\oplus \ker\chi$ the Riemannian metric $g=\chi\otimes\chi + \pi_b^*g_B$, being $g_B$ any metric on $B$ and $\pi_b$ the restriction of $\pi_*\colon TM\to TB$ to $\ker\chi$, which is an isomorphism by degree reasons. Then $\F$ is an isometric flow on $M$ whose Euler form is $d\chi=\pi^*\omega$. We shall use this construction in the following examples. First, we show that Theorem \ref{th:equivalencia} applies to new cases outside the contact category:

\begin{example}\label{ex:cp2cp2}
Let $B=\CP^2\sharp\CP^2$. Recall that its cohomology is given by: 
$$
\begin{aligned}
&\HB{0}=\R 	&\HB{2}&=<[a],[b]>=\R\oplus\R \\
&\HB{1}=\HB{3}=0 &\HB{4}&=<[a]^2>=<[b]^2>=\R,
\end{aligned}
$$
where $a$ and $b$ can be assumed to be integral. Choose $e=a$ (we could also take $e=b$).
Using Kobayashi's construction, we construct a 5-manifold $M$ with an isometric flow $\F$ whose  orbit space is $B$ and with Euler class $[e]=[a]\in\HB{2}$. As $B$ does not admit almost-complex structures,  (\cite[Prop.~1.3.1]{audin}), then it does not admit a symplectic structure, and $M$ cannot be a contact flow. On the other hand, as the Lefschetz maps:
$$
\xymatrix{
	L^2\colon H^0(B)\ar[r]^(0.6){\wedge [a]^2}& H^4(B) 
	&L^1\colon H^1(B)\ar[r]^(0.55){\wedge[a]} & H^3(B)
	&L^0\colon H^2(B)\ar[r]^(0.55){Id} & H^2(B)
}
$$
are isomorphisms, the flow satisfies $\THL{}$ and thus $\HL{}$.
So, $\F$ is a Lefschetz isometric flow that cannot be generated by the Reeb vector field of any contact structure.
\end{example}

\begin{remark}
	
	The manifold $\CP^2\sharp\CP^2$ is an example of a compact c-symplectic manifold (cohomologically symplectic manifold) which is not symplectic, as defined by Lupton and Oprea in \cite{LP2}. The authors also suggest a method of construction of other compact c-symplectic manifolds in \cite{LP1}. Therefore using such examples we can construct more Lefschetz flows which are not transversely symplectic, and thus not contact.
\end{remark}

\begin{example}
	It is very easy to construct isometric flows that are not Lefschetz flows. For example, let the circle $\ese^1$ act by multiplication on the first factor of $M=\ese^1\times B$ for any $B$ and call $\F$ the flow associated to that free action. Let $\pi_1$ and $\pi_2$ stand for the projections of $M$ onto $\ese^1$ and $B$, respectively.
	Consider the metric $g=\pi_1^*g_{\ese^1}+\pi_2^*g_B$ induced by any invariant metric in $\ese^1$ and any metric in $B$. Recall the characteristic form $\chi=i_Xg$, being $X$ the  unit vector field. It is straightforward to check locally that $d\chi=0$, i.e. the Euler form (hence, the Euler class) vanishes. So, $\F$ is not a Lefschetz flow.
	% Moreover, if $B=\CP^2\sharp\CP^2$, then we make sure that the example does not lie in the contact category for any metric. 
	%So we have both Lefschetz and non-Lefschetz flows outside the category of contact flows.
\end{example}
The following two lemmata will be useful to find integral closed forms within a family.

\begin{lemma}\label{lem:lattice}
	Let $n\in\N$. The only polynomial with real coefficients in $n$ variables that vanishes in a lattice of $\R^n$ is the zero polynomial.
\end{lemma}
\begin{proof}
	Let $q$ be a real polynomial vanishing in a lattice $\Lambda$ of $\R^n$. For any nonzero $z\in\Lambda$,  the restriction of $q$ to the line joining the origin and $z$ is a polynomial in one variable with an infinite amount of zeroes, hence zero. In particular, $q$ vanishes in any point of $\R^n$ with rational coefficients with respect to any basis of $\Lambda$, and by density, in all $\R^n$. 
\end{proof}

\begin{lemma}\label{lem:integral}
	Let $\HM{2}=<[\omega_1],\dots,[\omega_n]>$ and let $p$ be a nontrivial polynomial with real coefficients in $n$ variables. Then, the following subset contains an integral closed form:
	$$
	W_p=\Set{c_1\omega_1+\cdots+c_n\omega_n | c_i\in\R \text{ and } p(c_1,\dots,c_n)\ne 0}\subset\Omega^2(M).
	$$	
\end{lemma}
\begin{proof}
	Denote by  $H_b^2(M,\Z)$ the Betti part of $H^2(M,\R)$, i.e., the natural inclusion of $H^2(M,\Z)$ in $H^2(M,\R)$. Notice that $H^2_b(M,\Z)\cong \HM{2}$. Consider $\beta=\Set{[\alpha_1],\dots,[\alpha_n]}$ the image of a basis of $H^2_b(M,\Z)$ via the de Rham isomorphism. So, $\alpha_1,\dots,\alpha_n$ are integral cycles. As $\beta$ is a basis of $\HM{2}$, there exist $a_{ij}\in\R$ and forms $\gamma_i\in\Omega^1(M)$ such that
	\begin{equation}\label{eq:basis}
	\sum_{j=1}^n a_{ij}\omega_j=\alpha_i + d\gamma_i,
	\quad
	\text{for all }i\in\set{1,\dots,n}.
	\end{equation}
	Consider the polynomial in $n$ variables
	\begin{equation}\label{eq:poly}
	q(x_1,\dots,x_n)=p\left(\sum_{i=1}^n  a_{i1}x_i,\dots,\sum_{i=1}^n a_{in}x_i\right).
	\end{equation}
	By Lemma \ref{lem:lattice}, there exist integers $z_1,\dots,z_n$ such that $q(z_1,\dots,z_n)\ne 0$. Now we define
	%$$ 
	%c_j=\sum_{i=1}^n a_{ij}z_i,\quad 1\le j\le n
	%\quad
	%\text{and}
	%\quad
	%\omega=\sum_{i=1}^n z_i\sum_{j=1}^na_{ij}\omega_j=\sum_{j=1}^n c_j\omega_j
	%$$
	$$ 
	\omega=\sum_{i=1}^n z_i\sum_{j=1}^na_{ij}\omega_j=
	\sum_{j=1}^n\sum_{i=1}^n z_i a_{ij} \omega_j
	=\sum_{j=1}^n c_j\omega_j
	$$
	where $c_j=\sum_{i=1}^n a_{ij}z_i$ for  $j\in\set{1,\dots,n}$. On one hand, by \eqref{eq:poly}, we have $p(c_1,\dots,c_n)=q(z_1,\dots,z_n)\ne0$. Thus, $\omega\in W_p$. On the other hand, by \eqref{eq:basis}, we have
	$$
	\omega=\sum_{i=1}^n z_i\sum_{j=1}^na_{ij}\omega_j=\sum_{i=1}^n z_i\alpha_i + d\left(\sum_{i=1}^n z_i\gamma_i\right)
	$$
	and therefore, $\omega$ is an integral cycle.
\end{proof}

\begin{example}\label{ex:emeseis}
	In \cite{mino2} the authors construct K-contact flows which are not Lefschetz, and thus, not Sasakian. We can construct another example of the same kind (K-contact and not Lefschetz) as follows: consider $B$ the 6-dimensional solvmanifold defined in \cite{emeseis}. The authors construct a family of closed forms:
	$$
	\omega=a\omega_1 + b\omega_2 + c\omega_3,
	$$
	being $[\omega_i]$ certain generators of $\HB{2}$. They show that if $ac\ne 0$, then $\omega$ is a symplectic form, but $B$ fails to satisfy the Hard Lefschetz property. By applying Lemma \ref{lem:integral} to the polynomial $p(x,y)=xy$, we can suppose $\omega$ to be a closed integral form. So, the Kobayashi construction yields a K-contact manifold which is not Lefschetz.
\end{example}

\begin{example}
Let $B$ be one of the manifolds constructed in \cite{rendiconti}, which satisfy the following properties:
\begin{enumerate}[(i)]
	\item $B$ is a compact 6-solvmanifold of completely solvable type;
	\item $B$ admits a nondefinite K\"ahler metric, and thus a symplectic structure;
	\item $B$ satisfies the strong Lefschetz property;
	\item The Poincaré polynomial of $B$ is $P_B(t)=1+2t+5t^2+8t^3+5t^4+2t^5+t^6$.
\end{enumerate}

By (ii),  $B$ is a symplectic manifold. More precisely, in \cite[ p.~63]{rendiconti} it is shown that any linear combination $\omega$ of certain closed forms generating $\HB{2}$ with parameters $r, s, t, u, v$ is a valid symplectic form if the condition $r(sv+tu)\ne0$ is satisfied. By Lemma \ref{lem:integral}, we can choose $\omega$ to be integral and consider $M$ its associated Kobayashi construction. So, $M$ is a K-contact flow, and by (iii) it is a Lefschetz flow.

In \cite[Th.~3.3]{marisajoseba}, it is proved that $B$ does not admit a K\"ahler structure. We show here a shorter proof of this fact using a later result: by the Benson-Gordon's Theorem, proved by Hasegawa in \cite{hasegawa}, a manifold satisfying (i) admits a K\"ahler structure if and only if it is a complex torus. So, by (iv), $B$ does not admit a K\"ahler structure. In particular, $M$ does not admit a Sasakian structure.

So, we have a K-contact Lefschetz flow which does not admit a Sasakian structure. 
Other examples of isometric flows which are Lefschetz but not transversely K\"ahler 
can be constructed using examples of manifolds from  \cite{McD} which are symplectic and Lefschetz but not K{\"a}hler. For another example in dimension 5, see \cite{mino3}. 

\end{example}
We finish this paper by showing that the leftmost region of Fig.~\ref{fig:dibujico} is not empty.
\begin{example}\label{ex:solenoid}
As described in \cite[Remarque~3]{ghys-solenoid} the geodesic flow of the flat torus $\T^2$ induces a contact flow on its  unit tangent bundle whose basic cohomology has infinite dimension (and so, it cannot be isometric because the basic cohomology of a Riemannian foliation is finite-dimensional). More explicitly, it can be described as the usual contact flow on $\T^3=\R^3/\Z^3$ with 
%Reeb vector field  $X=\cos{t}\dd{x} + \sin{t}\dd{y}$ and 
contact form $\eta=\cos{2\pi t}\ dx + \sin{2\pi t}\ dy$. We have a fibration $\pi\colon M\to\R/\Z=\ese^1$ such that $\pi^{-1}([t])$ is foliated diffeomorphic to $\T^2$ with a linear flow of slope $t$, and is, thus, foliated by circles if $t\in\Q$ or by dense non-compact orbits otherwise. 
%This local geometry is not that of a Riemannian flow, so $\F$ is not isometric. 
It is easily checked that all basic 1-forms are written locally as $f(t)\ dt$. So, $H_B^1\cong H^1(\ese^1)=\R$ and every basic 1-form is closed, which yields $H^2_B=\Omega_B^2$ by degree reasons. It is straightforward to show that $\Omega^2_B\cong C^{\infty}([0,1],\{0,\frac{1}{2},1\})$ and hence, it has infinite dimension and the Lefschetz map between extremal degrees $L\colon H^0_B\to H^2_B$ is not an isomorphism, while by Remark~\ref{rem:extremal}, $\mathcal{R}^0$ is the graph of the isomorphism $H^0_M\cong H^3_M$. In contrast, the remaining Lefschetz map $L^0\colon H^1_B\to H^1_B$ is the identity isomorphism, but the relation $\mathcal{R}^1$ is not the graph of $H^1_M\cong H^2_M$, because in this example every closed basic 1-form is a primitive form by degree reasons and $PH_B^1=H^1_B$ is not isomorphic to $H^1_M$.  So, we have a contact flow which is not K-contact and fails to satisfy any of the Lefschetz properties, these failures happening at non-corresponding degrees. Notice that for isometric flows, the properties  $\THL{}$ and $\HL{}$  are intimately related at corresponding degrees.
\end{example}

\section*{Conflict of interest}
The authors declare that they have no conflict of interest.

\bibliography{nirebib}

\bibliographystyle{nireamsalpha}

\end{document}
